\documentclass{article}
\usepackage{graphicx} 
\usepackage{Test_1}

\title{From Navier–Stokes to Measure-Valued Euler–Brinkman Equations: Homogenization in the Inviscid Limit}
\author{Richard Höfer, Eleni Hübner-Rosenau}
\date{\today}
\usepackage{amsmath} 
\usepackage[maxbibnames=99,style=numeric,sortcites=true]{biblatex}
\usepackage{hyperref}
\begin{document}

\maketitle
\begin{abstract}
	    \noindent We study the behaviour of the solution $u_{\eps}$ to the Navier-Stokes equations with vanishing viscosity and a no-slip condition in a periodically perforated domain of grid size $\eps$. The behaviour depends on the particle size $\eps^\alpha$ and the viscosity $\eps^\gamma$ of the fluid (with $\alpha>1,\gamma>0$). We prove convergence to a measure-valued solution, provided that the local Reynolds number is small, in the subcritical ($\alpha+\gamma>3$) and critical ($\alpha+\gamma=3$) regime. In the first case, $u$ solves the Euler equations, whereas in the second case $u$ solves the Euler-Brinkman equations.
        This extends the results of \url{https://doi.org/10.1088/1361-6544/acfe56} for less regular limit functions. In particular, we also prove existence of measure-valued solutions to the Euler-Brinkman equations, along with weak-strong uniqueness.
	\end{abstract}
    \section{Introduction} Mixtures of fluids with many small particles, called \textit{suspensions}, appear in a lot in nature and industrial applications -- be it in the form of blood, dust storms or air-fuel injections in engines. Since tracking every particle is not feasible, one wants to find and study an \textit{effective model} for these mixtures (as has been done in for example \cite{Allaire.90a,FeiNamNec.16,LuYang.23,Mikelic91,MarKhr64,Tartar.80}). Moreover, many of the physically interesting fluids have low viscosity, so it also interesting so study the homogenization setting together with a vanishing viscosity coefficient. In this work, we will only consider stationary particles with a no-slip boundary condition.
    
    Some mathematical results dealing with such a combined homogenization and inviscid limit are \cite{Mikelic91, Mikelic95, MARUSICPALOKA200097}, where the problem is studied in a periodic geometry and the particle size is of the same order as the inter-particle distance, resulting in Darcy-type equations. Results dealing with smaller particles are \cite{lacave2015, Richard.23,höfer2026homogenizationnavierstokesequationsrandomly}: in \cite{lacave2015}, the subcritical $2$-dimensional regime is studied, in \cite{Richard.23}, all corresponding three regimes in the $3$-dimensional case are analyzed leading to the Euler, the Euler-Brinkman and Darcy's law in the subcritical, critical and supercritical case, respectively. In \cite{höfer2026homogenizationnavierstokesequationsrandomly}, the results of \cite{Richard.23} are extended to random particle configurations in (parts of) the critical and subcritical regime. \\
    The results \cite{lacave2015, Richard.23,höfer2026homogenizationnavierstokesequationsrandomly} are restricted to the cases where the limit function $u$ is sufficiently smooth (which is not necessary if the viscosity does not vanish, like in \cite{Allaire.90a,FeiNamNec.16,LuYang.23}). However, in light of the recent blow-up results for the Euler equations, smooth solutions can only be expected to exist for small times, even with smooth initial data. In this work, we extend the results of \cite{Richard.23} with less assumptions on the regularity of $u$, where we of course have to relax the notion of solution to the Euler(-Brinkman) equations. Namely, we prove that without assuming additional regularity of the limiting velocity field, up to a subsequence, the solutions converge to a measure-valued solution of the Euler--Brinkman equations. In particular, this provides existence for measure-valued solutions to the Euler-Brinkman system arising from a microscopic fluid model. 

    Requiring less regularity of the limit function is a natural question since the regime of low viscosity is generally characterized by turbulent flow (that is, flow with little regularity). Understanding  the regularity of fluids, even with initially smooth data, is in itself a very interesting question and has attracted much research. To this end, various notions of weaker solutions to the Euler equations have been introduced - for example the dissipative solutions already by Lions in the 90s \cite{Lions_1996}, and since then more general measure-valued solutions, dissipative weak solutions (confusingly, these are not the same as the ones from Lions) or energy-variational solutions. All of these concept have the property of global existence and weak-strong uniqueness. In fact, for the Euler equations, the latter three concepts have been shown to be equivalent, see \cite{eiter2026equivalencegeneralizedsolutionconcepts} (and the dissipative solutions of Lions are a special case of them). Since we want to study vanishing viscosity limits, we have decided to work with measure-valued solutions, where existence is usually shown by considering a vanishing viscosity sequence of the Navier-Stokes equations (see \cite{wiedemann2017weakstronguniquenessfluiddynamics,Gallenm_ller_2023,eiter2026equivalencegeneralizedsolutionconcepts} and references therein for further reading on this topic).\\
    This article is structured as follows: below, we first introduce the mathematical setting and state our main results. Then, we recall in Section \ref{sec:Young measures and mvs} important results for Young measures and measure-valued solutions to the Euler equations, as well as introducing measure-valued solutions to the Euler-Brinkman equations. In Section \ref{sec:corrector}, we recall some preliminary results from \cite{Richard.23} which we will use in Section \ref{sec:proof main results} to prove our main results.
\subsection{Setting and notation}
Let $\mathcal{T}\Subset B_{1/4}(0)$, the reference particle, be a closed smooth set such that $B_1(0)\setminus \mathcal{T}$ is connected and $0\in \mathring{\mathcal{T}}$. We then perforate the entire space $\R^3$ with holes on a periodic grid of size $\eps$, where we rescale the holes to size $\eps^{\alpha}$. To be precise, we define the perforated domain as
	\begin{align*}
		\Omega_\eps=\R^3\setminus \bigcup_{i\in\N} \mathcal{T}_i^\eps, \quad\qquad\mathcal{T}_i^\eps=x_i^\eps+\eps^\alpha \mathcal{T}.
	\end{align*}
	Here, $\alpha>1$ is a positive parameter and the centres $x_i^\eps$ are on the grid $\eps\Z^3$.\\
    \noindent We will study the behaviour of $u_\eps$, the solution to the incompressible Navier-Stokes equations with vanishing viscosity $\eps^\gamma$ (where $\gamma>0$), i.e.
	\begin{equation}\label{eq: pde for u_eps}
			\begin{aligned}
			\del_t u_\eps+(u_\eps\cdot \nabla)u_\eps-\eps^\gamma\Delta u_\eps+\nabla p_\eps&=0 &&\text{ in }[0,T]\times\Omega_\eps\\
			\divg(u_\eps)&=0&&\text{ in }[0,T]\times\Omega_\eps\\
			u_\eps|_{\partial\Omega_\eps}&=0&&\text{ on }[0,T]\\
			u_\eps(0,\cdot)&=u_{0,\eps}&&\text{ in }\Omega_\eps,
		\end{aligned}
	\end{equation}
	where $u_{0,\eps}\in L^2_{\sigma}(\Omega_\eps;\R^3)$ are given initial data. Here, $L^2_{\sigma}(\Omega_\eps;\R^3)$ denotes the space
    \begin{align*}
        L^2_{\sigma}(\Omega_\eps;\R^3)=\{v\in L^2(\Omega_\eps;\R^3)\:|\:\divg(v)=0\text{ (in the distributional sense)}, v\cdot n=0\text{ on }\partial\Omega_\eps\}
    \end{align*}
	
	It is known (see \cite{Sohr_2014} or \cite{Robinson_Rodrigo_Sadowski_2016} for an overview of this and other classical results) that a Leray solution exists, i.e. a weak solution satisfying the energy estimate (for almost all $t\in[0,T]$ with $t>s$, for almost all $s\geq 0$ including $s=0$)
	\begin{align}\label{eq:energy ineq}
		\frac{1}{2}\norm{u_\eps(t)}^2_{L^2(\Omega_\eps)}+\eps^\gamma\norm{\nabla u_\eps}^2_{L^2((s,t)\times\Omega_\eps)}\leq 	\frac{1}{2}\norm{u_{\eps}(s)}^2_{L^2(\Omega_\eps)}.
	\end{align} 
	We will show that when $\alpha+\gamma=3$ (critical scaling) and assuming that $u_{0,\eps}$ converges strongly to some limit function $u_0$, then $u_\eps$ converges to a function $u$ being a measure-valued solution (which we will define below) of the Euler-Brinkman equations 
	\begin{equation}\label{eq: pde for u}
		\begin{aligned}
					\del_t u+(u\cdot\nabla)u+\nabla p+\mathcal{R}u&=0\quad\text{in }[0,T]\times\R^3\\
			\divg(u)&=0\quad\text{in }[0,T]\times\R^3\\
			u(0,\cdot)&=u_0,\quad\text{in }\R^3.
		\end{aligned}
	\end{equation}
	 This extra linear term $\mathcal{R}u$, which represents the average effect that all the particles have on the fluid, is called the \textit{Brinkman force}. Here, $\mathcal{R}$ is the constant, symmetric, nonnegative resistance matrix, which can be determined by solving the Stokes equations with suitable boundary conditions on $\R^3\setminus\mathcal{T}$: let $w_k,q_k\in \dot{H}^1(\R^3)\times L^2(\R^3)/\R$ be the solutions to the exterior Stokes problem
	\begin{equation}
		\begin{aligned}
			-\Delta w_k+\nabla q_k&=0 \quad \text{in  }\R^3\setminus\mathcal{T},\\
			\divg(w_k)&=0\quad \text{in }\R^3\setminus\mathcal{T},\\
			w_k&=e_k\quad \text{on }\partial\mathcal{T},
		\end{aligned}
	\end{equation}
	where $e_k$ is the $k$-th unit vector in $\R^3$. Then we set
    \begin{align}\label{eq:resistance matrix}
        R_{jk}=\int_{\R^3\setminus\mathcal{T}}\nabla w_k:\nabla w_j\dif{x}.
    \end{align}
    In the subcritical case ($\alpha+\gamma>3$), we show that $u_\eps$ converges to a measure-valued solution of the Euler equations (corresponding to $\mathcal{R}=0$):
    	\begin{equation}\label{eq: pde for u subcrit}
		\begin{aligned}
					\del_t u+(u\cdot\nabla)u+\nabla p&=0\quad\text{in }[0,T]\times\R^3\\
			\divg(u)&=0\quad\text{in }[0,T]\times\R^3\\
			u(0,\cdot)&=u_0\quad\text{in }\R^3,
		\end{aligned}
	\end{equation}
\subsection{Statement of the main results}
\begin{thm}\label{thm:crit case}
    Let $u_\eps$ be a Leray solution to \eqref{eq: pde for u_eps} with $\alpha+\gamma=3$ and $\alpha>\gamma>0$. Then, assuming that $u_{0,\eps}\rightarrow u_0$ strongly in $L^2(\R^3)$, $u_\eps$ converges (up to a subsequence that we still label $u_\eps$) in the sense of Young measures to some $\nu\in L^\infty_w([0,T]\times\R^3;\mathcal{P}(\R^3))$, for which there exist $m,\mathcal{D}$ such that $(\nu, m, \mathcal{D})$ is a measure-valued solution to the Euler-Brinkman equations \eqref{eq: pde for u} (with $\mathcal{R}$ defined by \eqref{eq:resistance matrix}) with initial datum $u_0$.
\end{thm}
\begin{cor}\label{cor:existence Euler-Brinkman}
    Let $u_0\in L^2(\R^3)$. Then there exists a measure-valued solution $(\nu,m,\mathcal{D})$ to the Euler-Brinkman equations with initial datum $u_0$.
\end{cor}
\begin{thm}\label{thm:subcrit case}
    Let $u_\eps$ be a Leray solution to \eqref{eq: pde for u_eps} with $\alpha+\gamma>3$ with $\alpha>\gamma>0$. Then, assuming that $u_{0,\eps}\rightarrow u_0$ strongly in $L^2(\R^3)$, $u_\eps$ converges (up to a subsequence that we still label $u_\eps$) in the sense of Young measures to some $\nu\in L^\infty_w([0,T]\times\R^3;\mathcal{P}(\R^3))$, for which there exist $m,\mathcal{D}$ such that $(\nu, m, \mathcal{D})$ is a measure-valued solution to the Euler equations \eqref{eq: pde for u subcrit} with initial datum $u_0$.
\end{thm}
The strategy to prove these results is similar to the one used to prove existence of measure-valued solutions in for example \cite{wiedemann2017weakstronguniquenessfluiddynamics}. However, since in general test functions in $C_c^\infty(\R^3)$ are do not satisfy the no-slip boundary condition at the particles, we use the matrix-valued cutoff function from \cite{Richard.23}. Since we do not get quantitative convergence rates, the proof is in some aspects less technical than in \cite{Richard.23}.

\section{Young measures and measure-valued solutions}\label{sec:Young measures and mvs}
\subsection{Young measures and their properties}
For completeness, we briefly introduce the notion of Young measures and their properties, following a combination of the presentations in \cite{Gallenm_ller_2023} and \cite{wiedemann2017weakstronguniquenessfluiddynamics}.
\begin{defi}
    A \textit{Young measure} $\nu$ on $\R^d$ is an element of $L^\infty_w([0,T]\times\R^d;\mathcal{P}(\R^d))$ (denoting by the subscript $w$ the weak-$*$ topology). That is, it is a collection of probability measures $(\nu_{t,x})_{t\in[0,T],x\in\R^d}$ so that for every $g\in C_b(\R^d)$, the map
    \begin{align*}
        (t,x)\mapsto \int_{\R^d} g(\xi)\dif{\nu_{t,x}}(\xi)=:\langle g,\nu_{t,x}\rangle
    \end{align*}
    is measurable.
\end{defi}
We will write in the following for a (Young) measure $\nu$ and a measurable function $f:\xi\mapsto f(\xi)$
\begin{align*}
    \langle \nu, f(\xi)\rangle
\end{align*}
instead of $\langle \nu,f\rangle$ (such as $\langle \nu, \xi\rangle$ instead of $\langle \nu, \mathrm{id}\rangle$).

Young measures arise naturally as limits of $L^p$-bounded sequences (see below), so we introduce also the following useful notion of convergence, which behaves well under nonlinear compositions.
\begin{defi}
    A sequence of measurable functions $u_n$ is said to \textit{generate} a Young measure $\nu$ if for all $f\in C (\R^d)$ such that $f(u_n)$ is equiintegrable
    \begin{align}\label{eq:generating younf measure convergence}
        f(u_n(\cdot))\overset{\ast}{\rightharpoonup} \langle \nu_{(\cdot)},f\rangle\quad\text{in }L^\infty((0,T)\times\R^d).
    \end{align}
\end{defi}
The following fundamental theorem (which follows from Banach-Alaoglu) allows us to pass to limits for $L^2$-bounded sequences.
\begin{thm}[Fundamental theorem of Young measure theory]\label{thm:fundamental theorem young measures}
    Let $u_n: [0,T]\times \R^d\rightarrow \R^d$ be a family of functions bounded in $L^2([0,T]\times \R^d)$. Then, up to a non-relabeled subsequence, $u_n$ generates a Young measure $\nu$, which also satisfies
    \begin{align*}
        \langle \nu, |\xi|^2\rangle <\infty.
    \end{align*}
\end{thm}
\subsection{Measure-valued solutions}
We first recall the definition of measure-valued solutions to the Euler equations. We will denote by $\mathcal{M}(X,Y)$ the space of (signed) measures on $X$ with values in $Y$. Moreover, $T>0$ is a fixed, finite time.
\begin{defi}
A triple $(\nu,m,\mathcal{D})$ is a \emph{dissipative measure-valued solution}
of the incompressible Euler equations \eqref{eq: pde for u subcrit} with initial data
$u_0\in L^2_\sigma(\R^d)$ if
\begin{enumerate}
\item[(i)]
$\nu \in L^\infty([0,T]\times\R^d;\mathcal P(\mathbb R^d))$
is a Young measure;

\item[(ii)]
the measure
\[
m \in \mathcal M([0,T]\times\R^d;\mathbb R^{d\times d}_{\mathrm{sym}})
\]
satisfies
\[
m(\dif{t}\,\dif{x})=m_t(\dif{x})\otimes \dif{t}
\]
for some family of measures
\[
(t\mapsto m_t)\in
L^\infty_w([0,T];\mathcal M(\R^d;\mathbb R^{d\times d}));
\]

\item[(iii)]
the dissipation defect $\mathcal{D}\in L^\infty([0,T])$ satisfies $\mathcal{D}\geq 0$ and
\[
|m_t|(\R^d)\leq C\,\mathcal{D}(t)
\]
for some $C>0$ and for a.e. $t\in[0,T]$, as well as $D(0)=0$;

\item[(iv)]
the momentum equation
\begin{align*}
&\int_s^t \int_{\R^d}
(
\langle \nu_{t,x}, \xi\rangle \cdot \partial_t\varphi
+
\langle \nu_{t,x}, \xi\otimes \xi\rangle
:
\nabla\varphi
)\dif{x}\dif{t}
+
\int_s^t \int_{\R^d}
\nabla\varphi : m_t\dif{t}
\\
&=
\int_{\R^d} \langle \nu_{t,x},\xi\rangle\cdot\varphi(t,\cdot)-\langle \nu_{s,x},\xi\rangle\cdot\varphi(s,\cdot)\,dx
\end{align*}
holds for almost all $t\in [0,T]$ with $t\geq s$, for almost every $s\geq 0$ including $s=0$, and all
\[
\varphi\in C_c^\infty([0,T)\times{\R^d};\mathbb R^d)
\]
(or, by density, for all $\varphi\in C^1([0,T]\times{\R^d};\mathbb R^d)$) with $\operatorname{div}\varphi=0$;

\item[(v)]
the divergence-free condition
\[
\int_{\R^d}
\langle \nu_{t,x}, \xi\rangle \cdot \nabla\psi\dif{x}
=0
\]
holds for a.e.\ $t\in[0,T]$ and all
\[
\psi\in C_c^\infty({\R^d});
\]
\item[(vi)] the initial condition is attained via
\begin{align*}
    \langle\nu_{0,x},\xi\rangle=u_0(x);
\end{align*}
\item[(vii)]
the energy inequality
\[
\frac12
\int_{\R^d}
\langle \nu_{t,x}, |\xi|^2\rangle\,dx
+
\mathcal{D}(t)
\le
\frac12
\int_{\R^d}
\langle \nu_{s,x}, |\xi|^2\rangle\,dx
+
\mathcal{D}(s)
\]
holds for a.e.\ $0\leq s<t\in[0,T]$, including $s=0$.
\end{enumerate}
(All integrals are assumed to exist)
\end{defi}
\begin{rem}
    If $\nu_{t,x}=\delta_{u(t,x)}$ for some function $u$ and $\mathcal{D}=0$, $m=0$, then this definition recovers weak solutions to the Euler equations.
\end{rem}
\begin{rem}
    The requirement $D(0)=0$ matters only for the energy inequality, since altering a $L^\infty$ function at a single point does not change its properties otherwise.
\end{rem}
\begin{rem}
    Note that there are multiple similar definitions of measure-valued solutions in the literature; we chose to take the one from \cite[Definition 1]{Gallenm_ller_2023}, but append it with the requirement of the momentum equation and energy inequality holding for two arbitrary times $0\leq s<t$ (not just comparing to $t=0$), so that the results from \cite{eiter2026equivalencegeneralizedsolutionconcepts} are easily applicable to our results. Note that the above definition of measure-valued solutions is in fact also equivalent to the one considered in \cite{eiter2026equivalencegeneralizedsolutionconcepts}, which seems to be more general at first: every such solution easily is a solution in our sense by taking $m_t=\langle\nu^\infty,\theta\otimes\theta\rangle\lambda$ and $\mathcal{D}(t)=\frac12 \lambda_t(\R^3)$. We note that the proof of equivalence to the other solution concepts considered in \cite{eiter2026equivalencegeneralizedsolutionconcepts} works the same if one starts from a solution in our sense with more or less only notational adjustments.
\end{rem}
The following results are well-known (see \cite[Theorem 2]{Gallenm_ller_2023} and \cite[Theorem 3.4]{wiedemann2017weakstronguniquenessfluiddynamics}). The proof is done by considering a vanishing viscosity sequence of solutions to the Navier-Stokes equations.
\begin{thm}[Existence of measure-valued solutions to Euler]
    Let $d\geq 2$ and $u_0\in L^2_\sigma(\R^d)$. Then there exists a measure-valued solution to the Euler equations \eqref{eq: pde for u subcrit} with initial datum $u_0$.
\end{thm}
\begin{thm}[Weak-strong uniqueness of measure-valued solutions to Euler]
    Let $u_0$ be an initial datum such that there exists a strong solution $U\in C^1([0,T]\times\R^d)$ to the Euler equations \eqref{eq: pde for u subcrit} with initial datum $u_0$. Moreover, let $(\nu,m,\mathcal{D})$ be a measure-valued solution to the Euler equations with initial datum $u_0$. Then $(\nu,m,\mathcal{D})=(\delta_U,0,0)$.
\end{thm}
We define measure-valued solutions to Euler-Brinkman equivalently, adding a linear zeroth order term in the momentum equation:
\begin{defi}[Measure-valued solutions of the Euler-Brinkman equations]
Let $\mathcal{R}$ be the resistance matrix defined in \eqref{eq:resistance matrix}. A triple $(\nu,m,\mathcal{D})$ is a \emph{dissipative measure-valued solution}
of the incompressible Euler-Brinkman equations \eqref{eq: pde for u} with initial data
$u_0 \in L^2_\sigma({\R^d})$ if
\begin{enumerate}
\item[(i)]
$\nu \in L^\infty([0,T]\times{\R^d};\mathcal P(\R^d))$
is a Young measure;

\item[(ii)]
the measure
\[
m \in \mathcal M([0,T]\times{\R^d};\R^{d\times d})
\]
satisfies
\[
m(\dif{t}\,\dif{x})=m_t(\dif{x})\otimes \dif{t}
\]
for some family of measures
\[
(t\mapsto m_t)\in
L^\infty_w([0,T];\mathcal M({\R^d};\R^{d\times d}_\mathrm{sym}));
\]

\item[(iii)]
the dissipation defect $\mathcal{D}\in L^\infty(0,T)$ satisfies $\mathcal{D}\geq 0$ and
\[
|m_t|({\R^d})\leq C\,\mathcal{D}(t)
\]
for some $C>0$ and for a.e. $t\in[0,T]$, as well as $D(0)=0$;

\item[(iv)]
the momentum equation
\begin{align*}
    &\int_s^t \int_{\R^d}
(
\langle \nu_{t,x}, \xi\rangle \cdot \partial_t\varphi
+
\langle \nu_{t,x}, \xi\otimes \xi\rangle
:
\nabla\varphi
)+\mathcal{R}\langle \nu_{t,x}, \xi\rangle\phi\dif{x}\dif{t}
+
\int_s^t \int_{\R^d}
\nabla\varphi : m_t\dif{t}
\\
&=
\int_\Omega \langle \nu_{t,x},\xi\rangle\cdot\varphi(t,\cdot)-\langle \nu_{s,x},\xi\rangle\cdot\varphi(s,\cdot)\,dx
\end{align*}
holds for almost all $t\in [0,T]$ with $t\geq s$, for almost every $s\geq 0$ including $s=0$, and all
\[
\varphi\in C_c^\infty([0,T]\times{\R^d};\mathbb R^d)
\]
(or, by density, for all $\varphi\in C^1([0,T)\times{\R^d};\mathbb R^d)$) with $\operatorname{div}\varphi=0$ (here, $\mathcal{R}$ is a positive semidefinite matrix);

\item[(v)]
the divergence-free condition
\[
\int_{\R^d}
\langle \nu_{t,x}, \xi\rangle \cdot \nabla\psi\dif{x}
=0
\]
holds for a.e.\ $t\in[0,T]$ and all
\[
\psi\in C_c^\infty({\R^d});
\]
\item[(vi)] the initial condition is attained via
\begin{align*}
    \langle\nu_{0,x},\xi\rangle=u_0(x);
\end{align*}
\item[(vii)]
the energy inequality
\[
\frac12
\int_{\R^d}
\langle \nu_{t,x}, |\xi|^2\rangle\,dx
+
\mathcal{D}(t)
\le
\frac12
\int_{\R^d}
\langle \nu_{s,x}, |\xi|^2\rangle\,dx
+
\mathcal{D}(s)
\]
holds for a.e.\ $0\leq s<t$.
\end{enumerate}
(All integrals are assumed to exist)
\end{defi}
Below in Section \ref{sec:proof main results}, we will in particular prove existence of measure-valued solutions in $\R^3$ by considering the vanishing viscosity limit in a perforated domain. We also note that we have weak-strong uniqueness -- the proof works almost exactly like for the Euler equations:
\begin{thm}[Weak-strong uniqueness of measure-valued solutions to Euler-Brinkman]
    Let $\mathcal{R}$ be the resistance matrix defined in \eqref{eq:resistance matrix}. Let $u_0\in L^2_\sigma(\R^3)$ be an initial datum such that there exists a strong solution $U\in C^1([0,T]\times\R^d)$ to the Euler-Brinkman equations with initial datum $u_0$. Moreover, let $(\nu,m,\mathcal{D})$ be a measure-valued solution to the Euler-Brinkman equations \eqref{eq: pde for u} with initial datum $u_0$. Then $(\nu,m,\mathcal{D})=(\delta_U,0,0)$.
\end{thm}
\begin{proof}
    We note that by multiplying \eqref{eq: pde for u} by $U$ (which is admissible since $U$ is regular enough), we have since $\mathcal{R}\geq 0$
    \begin{align}
        \frac12\del_t \norm{U(t)}_2^2=-\int_{\R^3} \mathcal{R}U(t) U(t)\dx\leq 0.
    \end{align}
    Integrating this over time and using the energy inequality for $U$, we have for the following relative energy (choosing $s=0$)
    \begin{align*}
        E_{rel}(t)&:=\frac{1}{2}\int_{\R^d}\langle\nu_{t,x},|\xi-U(t,x)|^2\rangle\dx+\mathcal{D}(t)\\
        &=\frac{1}{2}\int_{\R^d}\langle\nu_{t,x},|\xi|^2\rangle+|U(t,x)|^2-2\langle\nu_{t,x},\xi\rangle U(t,x)\dx+\mathcal{D}(t)\\
        &\leq \int_{\R^d}|u_0|^2\dx-\int_{\R^d}\langle\nu_{t,x},\xi\rangle U\dx\\
        &=-\int_0^t\int_{\R^d}\langle\nu_{s,x},\xi\rangle \del_t U+\mathcal{R} \langle\nu_{s,x},\xi\rangle U+\nabla U:\langle\nu_{s,x},\xi\otimes\xi\rangle\dx\dif{s}-\int_0^t\int_{\R^d}\nabla U:\dif{m_s}\dif{s}.
    \end{align*}
    In the last equality, we used the weak momentum equation of $\nu$ with $U$ as a test function. Then, using that $U$ is a strong solution of the Euler-Brinkman equations with the symmetric matrix $\mathcal{R}$, we have
    \begin{align*}
        E_{rel}(t)&\leq-\int_0^t\int_{\R^d}-\divg(U\otimes U)\langle\nu_{s,x},\xi\rangle+\nabla U:\langle\nu_{s,x},\xi\otimes\xi\rangle\dx\dif{s}-\int_0^t\int_{\R^d}\nabla U:\dif{m_s}\dif{s}\\
        &=\int_0^t\int_{\R^d}\langle\nu_{s,x},[(\xi-U)\cdot\nabla] U(\xi-U)\rangle\dx\dif{s}-\int_0^t\int_{\R^d}\nabla U:\dif{m_s}\dif{s}\\
        &\leq C\int_0^t\int_{\R^d}\norm{\nabla U(s)}_{L^\infty}E_{rel}(s)\dif{s}.
    \end{align*}
    In the last inequality we have used that both $U$ and $\langle\nu,\xi\rangle$ are divergence-free, so that
    \begin{align*}
        \langle\nu,((\xi-U)\cdot \nabla) U \cdot U\rangle=0,\quad \divg(U\otimes U)=(U\cdot\nabla) U.
    \end{align*}
    Then, from Grönwall's inequality, it follows that $E_{rel}(t)=0$ for all $t$. Since both terms in the relative energy are nonnegative, this implies that $\mathcal{D}=0$ (implying that also $m=0$ since $|m_t|(\R^d)\leq C\mathcal{D}(t)$) and that $\nu_{t,x}=\delta_{U(t,x)}$ almost everywhere.
\end{proof}
\section{Construction and properties of the corrector $w^\eps$}\label{sec:corrector}
The strategy to prove Theorems \ref{thm:crit case} and \ref{thm:subcrit case} is similar to the one used to prove existence of measure-valued solutions in for example \cite{wiedemann2017weakstronguniquenessfluiddynamics}. However we face the problem that generally $\phi\in C_c^\infty([0,T)\times\R^3)$ is not an admissible test function for \eqref{eq: pde for u_eps}. Therefore, as in \cite{Richard.23}, we define an approximation
\begin{align}
    \phi^\eps=w^\eps \phi-B_\eps(\phi),
\end{align}
where $w^\eps=0$ on $\partial\Omega_\eps$ and $w^\eps\rightarrow\id$, and the operator $B_\eps$, which is a Bogovskii type operator, satisfying $\divg(B_\eps(\phi))=\divg(w^\eps\phi)=w^\eps:\nabla\phi$. Then, $\phi^\eps$ is an admissible test function and we will show that $\phi^\eps\rightarrow\phi$ strongly in $L^2([0,T];L^2(\R^3))$. In order to get the Brinkman term in the limit however, $w^\eps$ should also contain information about the boundary layers near the holes. In this Section, we explain in more detail how $w^\eps$ is constructed and prove some useful properties.

Around each particle centered on $x_i^\eps$ we have the cube with side lengths $\eps$, that is, $Q_i^\eps=x_i^\eps+[-\eps/2, \eps/2]^3$. We split this cube in four parts, namely
	\begin{align*}
		Q_i^\eps&=\mathcal{T}_i^\eps\cup D_i^\eps\cup C_i^\eps\cup K_i^\eps,\\
		C_i^\eps&:=B_{\eta_{\eps}/4}(x_i)\setminus \mathcal{T}_i^\eps,\\
		D_i^\eps&:=B_{\eta_{\eps}/2}(x_i)\setminus B_{\eta_{\eps}/4}(x_i),\\
		K_i^\eps&:= Q_i^\eps\setminus B_{\eta_{\eps}/2}(x_i).		
	\end{align*}
	Here,  
    \begin{align}
        \eta_\eps\in [\eps^\alpha,\eps]
    \end{align}
    is a parameter that we will choose later.
    
    We then construct the corrector function as follows: for $k=1,2,3$, let $w_k,q_k\in \dot{H}^1(\R^3)\times L^2(\R^3)/\R$ be the solutions to the exterior Stokes problem
	\begin{equation}
		\begin{aligned}
			-\Delta w_k+\nabla q_k&=0 \quad \text{in  }\R^3\setminus\mathcal{T},\\
			\divg(w_k)&=0\quad \text{in }\R^3\setminus\mathcal{T},\\
			w_k&=e_k\quad \text{on }\partial\mathcal{T},
		\end{aligned}
	\end{equation}
	where $e_k$ is the $k$-th unit vector in $\R^3$. Then we set
	\begin{align*}
		w_k^\eps(x)=e_k-w_k\left(\frac{x-x_i^\eps}{\eps^\alpha}\right), \: q_k^\eps=-\eps^{-\alpha}q_k\left(\frac{x-x_i^\eps}{\eps^\alpha}\right)\quad\text{in }C_i^\eps,\\
		-\Delta w_k^\eps+\nabla q_k^\eps=0,\: \divg(q_k^\eps)=0\quad\text{in } D_i^\eps,\\
		w_k^\eps=e_k,\:q_k^\eps=0\quad\text{in }\left(\R^3\setminus\bigcup_i B_{\eta_{\eps}/2}(x_i^\eps)
        \right).
	\end{align*}
	Note that for the definition inside $D_i^\eps$, we have to complement the PDE by suitable inhomogeneous boundary conditions, given by the definitions on $C_i^\eps$ and $K_i^\eps$. Moreover, recall that we assumed that $\mathcal{T}\Subset B_{1/4}$, so that $T_i^\eps\Subset B_{\eta_\eps/2}(x_i^\eps)$.
    
    \noindent In the following, we will write $w^\eps$ and $q^\eps$ for the matrix- respectively vector-valued functions that have $w_k^\eps$ respectively $q_k^\eps$ as their columns respectively entries. Moreover, we will write $A\leqc B$ whenever $A\leq CB$ for some constant $C$ depending only on the reference particle, $T$ or some (finite) norm of the smooth test function.\\  
	In the following Lemmas, we summarize important properties of this corrector. The results are taken from \cite[Section 2]{Richard.23}. First, we note some results on how well $w^\eps$ approximates the identity in various contexts.
\begin{lemma}[{\cite[Lemma 2.1]{Richard.23}}] The functions $w^\eps, q^\eps$ satisfy
		\begin{enumerate}[label=(\roman*)]
			\item $w^\eps\in W^{1,\infty}_0(\Omega_\eps;\R^{3\times 3})$, $q^\eps\in L^\infty(\Omega_\eps;\R^3)$, $\divg(w^\eps)=0$ and
			\begin{align}\label{eq:corrector w1,infty estimate}
				\norm{w^\eps}_{L^\infty(\R^3)}+\eps^\alpha\left(\norm{\nabla w^\eps}_{L^\infty(\R^3)}+\norm{q^\eps}_{L^\infty(\R^3)}\right)\leqc 1
			\end{align}
			\item For all $1\leq p<3$ it holds that
			\begin{align}\label{eq:corrector-Id Lp estimate p<3}
				\norm{(w^\eps-\id)\phi}_{L^p(\R^3)}^p\leqc \eps^{\alpha p-3}\eta_{\eps}^{3-p}\norm{\phi}_{W^{2,p}}^p\quad\text{for all }\phi\in W^{2,p}(\R^3).
			\end{align}
			Additionally,
			\begin{align}\label{eq:corrector-Id L3 estimate}
				\norm{(\id-w^\eps)\phi}_{L^3(\R^3)}^3&\leqc\eps^{3\alpha-3}|\log(\eps)|\norm{\phi}_{W^{2,3}}^3\quad\text{for all }\phi\in W^{2,3}(\R^3),\\
				\label{eq:corrector gradient L2 estimate ganzraum}
				\norm{\nabla w^\eps\phi}_{L^2(\R^3)}^2+\norm{q^\eps\phi}_{L^2(\R^3)}^2&\leqc \eps^{\alpha-3}\norm{\phi}_{H^2}^2\quad\text{for all }\phi\in H^2(\R^3),\\\label{eq:corrector gradient L1 estimate ganzraum}
				\norm{|\nabla w^\eps|^{1/2}\phi}_{L^2(\R^3)}^2+\norm{|q^\eps|^{1/2}\phi}_{L^2(\R^3)}^2&\leqc \eps^{\alpha-3}\eta_{\eps}\norm{\phi}_{H^2}\quad\text{for all }\phi\in H^2(\R^3).
			\end{align}
			\item For all $\phi\in H^1_0(\Omega_\eps)$
			\begin{align}\label{eq:corrector gradient L1 estimate Omega_eps}
				\norm{\phi |\nabla w^\eps|^{1/2}}_{L^2(\R^3)}^2+\norm{\phi |q^\eps|^{1/2}}_{L^2(\R^3)}^2&\leqc \eta_{\eps}\norm{\nabla\phi}_{L^2(\R^3)}^2.
			\end{align}
		\end{enumerate}
	\end{lemma}
    This corrector function mimics the behaviour of $u_\eps$ near the holes, which in the end will in the analysis cancel the Brinkman term from $u$.
    \begin{lemma}[{\cite[Lemma 2.2]{Richard.23}}]\label{lemma:right-hand side of stokes for corrector}
		We can write
		\begin{align*}
			-\Delta w^\eps+\nabla q^\eps=\eps^{\alpha-3}M_\eps-\gamma_\eps
		\end{align*}
		for some $M_\eps,\gamma_\eps\in W^{-1,\infty}(\R^3;\R^{3\times 3})$, where $\langle \gamma_\eps,v\rangle=0$ for all $v\in H^1_0(\Omega_\eps;\R^{3\times 3})$ and, for all $\psi\in H^1(\R^3;\R^3)$ and $\phi\in H^3(\R^3;\R^3)$, we have
		\begin{equation}\label{eq:estimate for M_eps}
			\begin{aligned}
				|\langle \phi(M_\eps-\rho\mathcal{R}),\psi\rangle|\leqc \left(\eta_{\eps}^{-1/2}\eps^{3/2}\norm{\psi}_{H^1(\R^3)}+\eta_{\eps}^{-1}\eps^\alpha\norm{\psi}_{L^2(\R^3)}\right)\norm{\phi}_{H^3(\R^3)}
			\end{aligned}
		\end{equation}
	\end{lemma}
    In order to make our test function divergence-free, we construct a Bogovskii-type operator on the perforated domain $\Omega_\eps$.
    \begin{lemma}[{\cite[Lemma 2.3]{Richard.23}}]
		For all $1<p<\infty$, there exists a  linear operator $B_\eps: W^{1,p}(\R^3;\R^3)\rightarrow W^{1,p}_0(\Omega_\eps;\R^3)$ such that for all $\phi\in W^{1,p}(\R^3;\R^3)$  that satisfy $\divg(\phi)=0$ we have
        \begin{align*}
            \divg(B_\eps(\phi))=w^\eps:\nabla\phi
        \end{align*}
        and
		\begin{align}\label{eq:estimate Bogovskii gradient}
			\norm{\nabla B_\eps(\phi)}_{L^p(\R^3)}^p&\leqc \norm{(I-w^\eps):\nabla\phi}_{L^p(\R^3)}^p,\\\label{eq:estimate Bogovskii} \norm{B_\eps(\phi)}_{L^p(\R^3)}^p&\leqc\eta_\eps^p \norm{(I-w^\eps):\nabla\phi}_{L^p(\R^3)}^p.
		\end{align}
	\end{lemma}
    Finally, we note the following Poincaré inequality (which follows by summing the Poincaré inequality on each ball $B_{\eta_\eps}(x_i^\eps)$).
    \begin{lemma}
        Let $1\leq p<\infty$ and let $v\in W^{1,p}_0(\Omega_\eps)$ be a function whose support is contained in $\bigcup_i B_{\eta_\eps}(x_i^\eps)$. Then
        \begin{align}\label{eq:Poincare in cup bEta}
            \norm{v}_{L^p(\R^3)}\leqc\eta_\eps\norm{\nabla v}_{L^p(\R^3)}.
        \end{align}
    \end{lemma}
\section{Homogenization of measure-valued solutions}\label{sec:proof main results}
\begin{proof}[Proof of Corollary \ref{cor:existence Euler-Brinkman}]
    Use Theorem \ref{thm:crit case} with, for example, $u_{0,\eps}=w^\eps u_0+B_\eps(u_0)\rightarrow u_0$ (even strongly) in $L^2(\R^3)$.
\end{proof}
\begin{proof}[Proof of Theorems \ref{thm:crit case} and \ref{thm:subcrit case}]
    We will focus on the proof of Theorem \ref{thm:crit case} and indicate at the end which parts need to be adapted in the subcritical case. In this proof, we will write (for example) $L^p(L^q)$ instead of of $L^p((0,T);L^q(\R^3))$. Also, we will still denote by $u_\eps$ the extension of $u_\eps$ to $\R^3$ by defining it to be zero on the particles (and likewise for other functions defined a priori only on $\Omega_\eps$).\\
    First, we observe that convergence of $u_\eps$ to some Young measure $\nu$ follows from the $L^2(L^2)$-boundedness, implying that we can consider $f(\xi)=\xi$ in \eqref{eq:generating younf measure convergence} by the de la Vallée-Poussin criterion:
    \begin{lemma}[de la Vallée-Poussin criterion,{\cite[Theorem 6.19]{klenke2014probability}}]
        Let $(\Omega,\mathcal{A},\mu)$ be a finite measure space. Then a family $\mathcal{F}\subset L^1(\mu)$ is equiintegrable if and only if there exists a function $H:[0,\infty)\rightarrow[0,\infty)$ with $H(\xi)/\xi\rightarrow\infty$ and 
        \begin{align*}
            \sup_{f\in\mathcal{F}}\int H(|f|)\dif{\mu}<\infty.
        \end{align*}
    \end{lemma}
    Now, take $\phi\in C_c^\infty([0,T]\times\R^d;\R^d)$ with $\divg(\phi)=0$, and $\psi\in C_c^\infty([0,T]\times\R^d)$. Moreover, let $w^\eps$ be the usual matrix-valued corrector as defined above. Then we define
    \begin{align}
        \phi^\eps=w^\eps\phi+B_\eps(\phi).
    \end{align}
    Here, $B_\eps$ is the Bogovskii operator on $\Omega_\eps$ as introduced above. We will choose the parameter $\eta_\eps$ such that 
    \begin{align}
    \eps^{3-\gamma}\ll\eta_\eps\ll\eps^\gamma
    \end{align}
    (for example, $\eta_\eps=\eps^{3/2}$ would work). Below, we will still keep $\eta_\eps$ as a parameter though, to see where this requirement comes from.  
    
    \textbf{Step 1: Divergence-free condition. }Since $u_\eps$ is divergence-free and $u_\eps=0$ on $\del\Omega_\eps$, we have
    \begin{align*}
        0=\int_{\Omega_\eps} u_\eps\cdot\nabla \psi^\eps\dif{x}\dif{s}=\int_{\R^3} u_\eps\cdot\nabla \psi^\eps\dif{x}\dif{s}
    \end{align*}
    By convergence of $u_\eps$, we have
    \begin{align*}
        \int_{\R^3} u_\eps\cdot \nabla \psi\rightarrow \int_{\R^3} \langle \nu_{t,x}, \xi\rangle\cdot \nabla \psi.
    \end{align*}
    This shows that 
    \begin{align}
        \int_{\R^3} \langle \nu_{t,x}, \xi\rangle\cdot \nabla \psi=0.
    \end{align}

    \textbf{Step 2: Momentum equation}\\
    We fix $0\leq s<t$ and consider the limit in each of the terms in the momentum equation. 
    
    \textbf{Step 2.1: Time derivative.} For the first term, we compute the time derivative of $\phi^\eps$:
    \begin{align*}
        \del_t\phi^\eps=w^\eps\del_t\phi-B_\eps(\del_t\phi),
    \end{align*}
    so that, by using \eqref{eq:corrector-Id Lp estimate p<3}, implying also that $B_\eps(\del_t\phi)\rightarrow 0$ strongly in $L^\infty(L^2)$ by \eqref{eq:estimate Bogovskii}
    \begin{align*}
        &\int_s^t\int_{\R^3} u_\eps\del_t\phi^\eps-\langle \nu_{\tau,x},\xi\rangle\del_t\phi\dif{x}\dif{\tau}\\
        &=\int_s^t\int_{\R^3} (u_\eps-\langle \nu_{\tau,x},\xi\rangle)\del_t\phi\dif{x}\dif{s}+\int_s^t\int_{\R^3} u_\eps((w^\eps-I)\del_t\phi-B_\eps(\del_t\phi))\dif{x}\dif{\tau}\\
        &\leq o(1)+T\norm{u_\eps}_{L^\infty(L^2)}\left[\norm{(w^\eps-\id)\del_t\phi}_{L^\infty(L^2)}+\norm{(w^\eps-\id):\nabla\del_t\phi}_{L^\infty(L^2)}\right]\rightarrow 0.
    \end{align*}
    \textbf{Step 2.2: Diffusion term.} For the diffusion term, we have, using integration by parts along with Lemma \ref{lemma:right-hand side of stokes for corrector} (note that the pressure disappears since $\phi^\eps$ is divergence-free and $\langle\gamma_\eps,\phi^\eps\rangle=0 $ as $\phi^\eps\in H^1_0(\Omega_\eps)$),
    \begin{align*}
        &\eps^\gamma\int_s^t\int_{\R^3} \nabla u_\eps:\nabla \phi^\eps\dif{x}\dif{\tau}\\
        &=\eps^\gamma\int_s^t \eps^{\alpha-3}\langle M_\eps u_\eps,\phi\rangle+\int_{\R^3}-u_\eps\nabla w^\eps\nabla\phi+\nabla u_\eps:w^\eps\nabla\phi+\nabla u_\eps:\nabla B_\eps(\phi)\dif{x}\dif{\tau}\\
        &= I_1+I_2+I_3+I_4.
    \end{align*}
    We consider each term separately. For the first integral, we use that $\alpha+\gamma=3$ and that $\eta_\eps\geq \eps^\alpha$, combined with \eqref{eq:estimate for M_eps} and the energy inequality for $u_\eps$ (implying that $\norm{\nabla u_\eps}_{L^2(L^2)}\leq C\eps^{-\gamma/2}$ and $\norm{u_\eps}_{L^\infty(L^2)}\leq C$:
    \begin{align*}
        \left|I_1-\int_s^t\int_{\R^3} \mathcal{R}\langle\nu_{\tau,x},\xi\rangle\phi\dif{\tau}\right|&=\left|\int_s^t\int_{\R^3}\mathcal{R}(u_\eps-\langle\nu_{\tau,x},\xi\rangle)\phi\dif{x}+\langle(M_\eps-\mathcal{R})u_\eps,\phi\rangle\dif{\tau}\right|\\
        &\leq o(1)+C\eta_\eps^{-1/2}\eps^{3/2}\norm{\nabla u_\eps}_{L^2(L^2)}+\eta_\eps^{-1}\eps^\alpha\norm{u_\eps}_{L^\infty(L^2)}\\
        &\leq o(1)+C\eta_\eps^{-1/2}\eps^{(3-\gamma)/2},
    \end{align*}
    which vanishes since we choose $\eta_\eps\gg \eps^{3-\gamma}=\eps^\alpha$. Here, we also used that $\eta_\eps^{-1/2}\eps^{\alpha/2}\leq 1$ in the last inequality.\\
    For the second integral, we split up $\nabla w^\eps$ and apply \eqref{eq:corrector gradient L1 estimate Omega_eps} and \eqref{eq:corrector gradient L1 estimate ganzraum}. Then we use again that $\norm{\nabla u_\eps}_{L^2(L^2)}\leqc \eps^{-\gamma/2}$:
    \begin{align*}
        |I_2|\leq \eps^\gamma \norm{|\nabla w^\eps|^{1/2}u_\eps}_{L^2(L^2)}\norm{|\nabla w^\eps|^{1/2}\nabla\phi}_{L^2(L^2)}\leq C\eta_\eps^{1/2}\eps^{\gamma/2}\eps^{(\alpha-3)/2}\eta_\eps^{1/2}=C\eta_\eps\rightarrow 0.
    \end{align*}
    Similarly, 
    \begin{align*}
        |I_3|\leq \eps^\gamma\norm{w^\eps}_{L^\infty}\norm{\nabla u_\eps}_{L^2(L^2)}\norm{\nabla\phi}_{L^2(L^2)}\leq C\eps^{\gamma/2}\rightarrow 0,
    \end{align*}
    and, using additionally \eqref{eq:estimate Bogovskii gradient},
    \begin{align*}
        |I_4|\leq \eps^\gamma \norm{\nabla u_\eps}_{L^2(L^2)}\norm{\nabla B_\eps(\phi)}_{L^2(L^2)}\leq C\eps^{\gamma/2}\eta_\eps^{1/2}\eps^{\alpha-3/2}\rightarrow 0.
    \end{align*}
    \textbf{Step 2.3: Convective term.} For the convective term, we note that we can mostly use the standard approach: for $u_\eps\otimes u_\eps$ we only have an $L^\infty(L^1)$ bound, which implies that (up to another subsequence that we will still denote $u_\eps$)
    \begin{align}\label{eq:convergence of nonlinear term to measure}
        u_\eps\otimes u_\eps-\langle\nu,\xi\otimes\xi\rangle\overset{\ast}{\rightharpoonup} m
    \end{align}
    for some measure $m\in\mathcal{M}([0,T]\times\R^3)$. By the disintegration theorem and the $L^\infty(L^1)$ bound, we can write
    \begin{align*}
        m(\dif{x}\dif{t})=m_t(\dif{x})\otimes \dif{t}
    \end{align*}
    for some family of measures $m_t\in \mathcal{M}(\R^3)$ that is uniformly bounded.
    Therefore, we can write
    \begin{align*}
        &\int_s^t\int_{\R^3} u_\eps\otimes u_\eps:\nabla\phi^\eps\dx\dif{\tau}-\int_0^t\int_{\R^3} \langle\nu_{\tau,x},\xi\otimes\xi\rangle:\nabla\phi\dx+\nabla\phi:\dif(m_\tau)\dif{\tau}\\
        &=o(1)+\int_s^t\int_{\R^3} (w^\eps-\id) u_\eps\otimes u_\eps:\nabla\phi\dx\dif{\tau}+\int_s^t\int_{\R^3} u_\eps\otimes u_\eps:\nabla w^\eps\phi\dx\dif{\tau}\\
        &=o(1)+I_{conv}^1+I_{conv}^2.
    \end{align*}
    For the first integral, we use the Poincaré inequality \eqref{eq:Poincare in cup bEta} (noting that $w^\eps-\id$ is supported in $\bigcup_i B_{\eta_\eps}(x_i^\eps)$), along with \eqref{eq:corrector gradient L1 estimate Omega_eps}:
    \begin{align*}
        |I_{conv}^1|&\leq \norm{\nabla\phi}_{L^\infty(L^\infty)}\eta_\eps\norm{\nabla [(w^\eps-\id)u_\eps\otimes u_\eps]}_{L^1(L^1)}\\
        &\leqc \eta_\eps\left(\norm{|\nabla w^\eps|^{1/2}u_\eps}_{L^2(L^2)}^2+\norm{(w^\eps-\id)}_{L^\infty}\norm{u_\eps}_{L^\infty(L^2)}\norm{\nabla u_\eps}_{L^2(L^2)}\right)\\
        &\leqc \eta_\eps(\eta_\eps\norm{\nabla u_\eps}_{L^2(L^2)}^2+\norm{\nabla u_\eps}_{L^2(L^2)})\\
        &\leqc \eta_\eps(\eta_\eps\eps^{-\gamma}+\eps^{-\gamma/2}).
    \end{align*}
    For the second integral, we use again \eqref{eq:corrector gradient L1 estimate Omega_eps} to estimate
    \begin{align*}
        |I_{conv}^2|\leqc \eta_{\eps}\norm{\nabla u_\eps}_2^2.
    \end{align*}
    Therefore, since we choose $\eta_\eps\ll \eps^\gamma$, both of these integrals converge to zero. 

    \textbf{Step 2.4: Initial condition.} For almost every $t>0$, we may pass to the limit in
    \begin{align*}
        \int_{\R^3}\langle\nu_{t,x},\xi\rangle \phi(t)\dif{x}
    \end{align*}
    since $u_\eps\overset{\ast}{\rightharpoonup}\langle \nu,\xi\rangle$ in $L^\infty([0,T]\times \R^3)$. Moreover, assuming that the initial data $u_{0,\eps}\rightarrow u_0$ strongly in $L^2(\R^3)$, we also have that since $\phi^\eps(0)\rightarrow \phi(0)$ strongly in $L^2$,
    \begin{align*}
        \int_{\R^3} u_{0,\eps}\phi^\eps(0)\rightarrow \int_{\R^3} u_0\phi(0).
    \end{align*}
    So in particular, $s=0$ is allowed in the momentum equation. 
    
    Therefore, with our previous choice of $\eta_\eps$, we can pass to the limit in all terms, implying that $\langle\nu,\xi\rangle$ satisfies the momentum equations for Euler-Brinkman. 
    
    \textbf{Step 3: Energy inequality.} Finally, we note that the dissipation defect
    \begin{align}\label{eq:definition of D as trace}
        \mathcal{D}(t):=\frac{1}{2}\mathrm{tr}(m_t)(\R^d)
    \end{align}
    is positive and the energy inequality is satisfied by exactly the same calculations as in the standard case of the Euler equations. For completeness, we briefly give the argument.

    First, since $m_t$ is uniformly bounded as a measure, $\mathcal{D}\in L^\infty([0,T])$. Also, $\mathcal{D}$ is nonnegative since for every $\chi\geq 0$ and $M>0$
    \begin{align}
        &2\int_0^t\chi(s)\mathcal{D}(s)\dif{s}=\lim_{\eps\rightarrow 0} \int_0^t\int_{\R^3}\chi(s)\mathrm{tr}(u_\eps\otimes u_\eps)-\langle\nu,\xi\otimes\xi\rangle)\dx\dif{s}\label{eq: equality for trace of m and energy1}\\
        &=\lim_{\eps\rightarrow 0} \int_0^t\int_{\R^3}\chi(s)\mathrm{tr}[(u_\eps\otimes u_\eps)-\langle\nu,\xi\otimes\xi\rangle)]\dx\dif{s}\\
        &=\lim_{\eps\rightarrow 0} \int_0^t\int_{\R^3}\chi(s)[|u_\eps|^2-\langle\nu,|\xi|^2\rangle]\dx\dif{s}\label{eq: equality for trace of m and energy3}\\
        &=\lim_{\eps\rightarrow 0} \int_0^t\int_{\R^3}\chi(s)[\min\{M,|u_\eps|^2\}-\min\{M,\langle\nu,|\xi|^2\rangle\}]\dx\dif{s}\nonumber\\
        &\quad+\lim_{\eps\rightarrow 0} \int_0^t\int_{\R^3}\chi(s)[\max\{0,|u_\eps|^2-M\}-\max\{0,\langle\nu,|\xi|^2\rangle-M\}]\dx\dif{s}\nonumber.
    \end{align}
    Now, we can use the fundamental theorem with the bounded test function $g(y)=\min\{M,|y|^2\}$ to see that the first term converges to zero. Also, since $\max\{0,|u_\eps|^2-M\}\geq 0$, we have
    \begin{align*}
        &2\int_0^t\chi(s)\mathcal{D}(s)\dif{s}=\lim_{\eps\rightarrow 0} \int_0^t\int_{\R^3}\chi(s)\mathrm{tr}(u_\eps\otimes u_\eps)-\langle\nu,\xi\otimes\xi\rangle)\dx\dif{s}\\
        &\geq -\int_0^t\int_{\R^3}\chi(s)\max\{0,\langle\nu,|\xi|^2\rangle-M\}\dx\dif{s}.
    \end{align*}
    However, this last term becomes arbitrarily close to zero by dominated convergence for $M\rightarrow\infty$ (which is applicable since $\max\{0,|y|^2-M\}\leq |y|^2$, the latter of which is integrable with respect to $\nu$). Therefore, $\mathcal{D}\geq 0$. 

    Finally, we note that since $u_\eps$ is a Leray-Hopf solution the the Navier-Stokes equation, it holds that for almost all $0\leq s<t$, including $s=0$,
    \begin{align*}
        \int_{\R^3}|u_\eps(t)|^2\dx\leq\int_{\R^3}|u_\eps(s)|^2\dx.
    \end{align*}
    Then, we pass to the limit by using the above calculation (see equalities \eqref{eq: equality for trace of m and energy1}-\eqref{eq: equality for trace of m and energy3}) along with the definition of $\mathcal{D}$ via \eqref{eq:definition of D as trace}. Moreover, we note that for $s=0$, we get the desired energy inequality with $D(0)=0$ by the strong convergence in $L^2(\R^3)$ of the initial data.

    \textbf{Adaptations in the subcritical case. } If instead we have $\alpha+\gamma>3$, almost all terms remain unchanged, with the exception of one of the diffusion terms: indeed, now we have that
    \begin{align*}
        I_1&=\eps^{\alpha+\gamma-3}\int_s^t \langle M_\eps u_\eps,\phi\rangle\dif{\tau}=\eps^{\alpha+\gamma-3}\int_s^t \langle (M_\eps-\mathcal{R}) u_\eps,\phi\rangle\dif{\tau}+\eps^{\alpha+\gamma-3}\int_s^t \langle \mathcal{R} u_\eps,\phi\rangle\dif{\tau}\\
        &\leqc \eps^{\alpha+\gamma-3}\eta_\eps^{-1/2}\eps^{(3-\gamma)/2}+\eps^{\alpha+\gamma-3}\rightarrow 0.
    \end{align*}
\end{proof}
\section*{Acknowledgments}
The second author was funded by the Deutsche Forschungsgemeinschaft (DFG, German
Research Foundation) Research Training Group 2339 “Interfaces, Complex Structures, and
Singular Limits in Continuum Mechanics – Analysis and Numerics”. The support is gratefully acknowledged.\\
The authors also would like to thank Emil Wiedemann for his helpful discussions and feedback, and Franck Sueur for suggesting to investigate this topic to the first author.
	\printbibliography
\end{document}